\documentclass[12pt]{amsart}
\usepackage{amssymb,amsmath,amsfonts,pictex,graphicx,ifthen,subfigure,fullpage,etex}

\DeclareMathOperator{\id}{id}

\newcommand{\SH}{\stackrel{\textrm{SH}}{\simeq}}
\newcommand{\CE}{\stackrel{\textrm{CE}}{\simeq}}
\newcommand{\fCE}{\stackrel{\textrm{CE}^{\textrm{f}}}{\simeq}}
\newcommand{\HE}{\stackrel{\textrm{HE}}{\simeq}}

\newcommand{\eS}{\mathcal{S}}

\newcommand{\then}{\;\Rightarrow\;}

\newtheorem{theorem}{Theorem}

\newtheorem{prop}{Proposition}[section]

\newtheorem{facts}[prop]{Facts}

\newtheorem{defn}[prop]{Definition}

\newtheorem*{ack}{Acknowledgements}
\newtheorem*{claim}{Claim}

\newtheorem*{Thm}{Theorem}

\begin{document}

\title[CAT(0) Boundaries of Groups which Split]{All CAT(0) Boundaries of a Group of the form \(H\times K\) are CE Equivalent}
\begin{abstract}
M. Bestvina has shown that for any given torsion-free CAT(0) group $G$, all of its boundaries are
\textit{shape equivalent}.  He then posed the question of whether they satisfy
the stronger condition of being \textit{cell-like equivalent}.
In this article we prove that the answer is ``Yes'' in the situation where the group in
question splits as a direct product with infinite factors.
We accomplish this by proving an interesting theorem in shape theory.
\end{abstract}

\author{Christopher Mooney}
\address{Department of Mathematical Sciences, University of Wisconsin-Milwaukee, Milwaukee, Wisconsin 53201}
\email{cpmooney@uwm.edu}
\date{July 27, 2007}
\subjclass{57M07, 20F65, 54C56}
\maketitle

\section{Introduction}
\label{sec:intro}
The CAT(0) condition is a geometric notion of nonpositive curvature, similar to the
definition of Gromov $\delta$-hyperbolicity.  A geodesic space $X$ is called CAT(0)
if it has the property that geodesic triangles in $X$ are ``no fatter'' than geodesic
triangles in euclidean space (see \cite[Ch II.1]{BH} for a precise definition).
The \textit{visual} or \textit{ideal boundary} of $X$, denoted \(\partial X\),
is the collection of endpoints of geodesic rays
emanating from a chosen basepoint.
It is well-known that \(\partial X\) is well-defined and independent of choice of basepoint.
Furthermore, when given the cone topology, \(X\cup\partial X\) is a $Z$-set compactification for $X$.
A group $G$ is called CAT(0) if it acts geometrically
(i.e. properly discontinuously and cocompactly by isometries) on some CAT(0) space $X$.
In this setup, we call $X$ a CAT(0) $G$-space
and \(\partial X\) a CAT(0) boundary of $G$.
We say that a CAT(0) group $G$ is \textit{rigid} if it has only one
topologically distinct boundary.\\

It is well-known that if $G$ is negatively curved (acts geometrically on a Gromov
$\delta$-hyperbolic space) or if $G$ is free abelian then $G$ is rigid.  Apart from
this, little is known concerning rigidity of groups.  P.L. Bowers and K. Ruane
showed that if $G$ splits as the product of a negatively curved group with
a free abelian group, then $G$ is rigid (\cite{BR}).
Ruane proved later in \cite{Ru} that if $G$ splits as a product of two negatively curved groups, then $G$
is rigid.  T. Hosaka has extended this work to show that in fact it suffices to know
that $G$ splits as a product of rigid groups (\cite{Ho}).
Another condition which guarantees
rigidity is knowing that $G$ acts on a CAT(0) space with isolated flats, which was proven
by C. Hruska in \cite{Hr}.\\

Not all CAT(0) groups are rigid, however: C. Croke and B. Kleiner constructed in \cite{CK}
an example of a non-rigid CAT(0) group $G$.  Specifically, they showed that $G$ acts on two different
CAT(0) spaces whose boundaries admit no homeomorphism.
J. Wilson proved in \cite{Wi} that this same group has \textit{uncountably many} boundaries.
Furthermore, it is shown in \cite{Mo} that the knot group $G$ of any connected sum of two non-trivial
torus knots has uncountably many CAT(0) boundaries.  For a collection of non-rigid
CAT(0) groups with boundaries of higher dimension, see \cite{MoCK}.\\

On the other end of the spectrum, it has been proven by M. Bestvina in \cite{Be}
that for any torsion-free CAT(0) group, all of its boundaries are \textit{shape equivalent}.
He then posed the question of whether they satisfy
the stronger condition of being \textit{cell-like equivalent}.
Bestvina's question
has been answered in part by R. Ancel, C. Guilbault, and J. Wilson, who showed in \cite{AGW}
that all the currently known boundaries of Croke and Kleiner's original group satisfy
this property; they are all cell-like equivalent to the Hawaiian earring.\\

In this article, we give further evidence in favor of Bestvina's conjecture by proving
the following theorem.\\

\begin{theorem}
\label{thm:celllike}
Let $G$ be a CAT(0) group which splits as a product $H\times K$
where $H$ and $K$ are infinite.
Then all CAT(0) boundaries of $G$ are cell-like equivalent through finite dimensional compacta.
\end{theorem}
Contrasting this with Hosaka's result, no assumption needs to be made about the factor groups.\\

In order to prove Theorem \ref{thm:celllike}, we first prove an interesting result in shape theory.
In \cite{Has}, Hastings proves that if two spaces are shape equivalent, then their
suspensions are cell-like equivalent.
The proof of this next theorem was inspired by a geometric proof of Hastings' theorem shown to
the author by Craig Guilbault.
\begin{theorem}
\label{thm:SH2CE}
Joins of shape equivalent compacta are cell-like equivalent.
That is, if \(X\SH X'\) and \(Y\SH Y'\), then
\[
	X\ast Y\CE X'\ast Y'
\]
Furthermore, if these four compacta are finite dimensional, then the cell-like equivalence can be realized
through finite dimensions.
\end{theorem}
Here $\ast$ denotes the join operation, $\SH$ denotes shape equivalence, and $\CE$
denotes cell-like equivalence.  For us, the term ``compactum'' means a compact metric space.\\

\begin{ack}
The work contained in this paper will be published as one part of the author's Ph.D. thesis written
under the direction of Craig Guilbault at the University of Wisconsin-Milwaukee.  The author would also
like to thank Ric Ancel, Chris Hruska, Boris Okun, and Tim Schroeder for helpful conversations.
\end{ack}

\label{sec:cellike}
\section{Equivalence of Compacta}
\label{sec:equiv}

\subsection{Shape Equivalence}
A number of definitions of shape equivalence have been given (see \cite{Bo} and \cite{MS}).
We will use the following equivalent definition, which is due to Chapman (\cite[Sec VI]{Ch}).\\

\begin{defn}
\label{defn:shape}
We say that two compacta $X$ and $Y$ are \textit{shape equivalent} and write \(X\SH Y\)
if when $X$ and $Y$ are imbedded as Z-sets in the Hilbert cube \(Q\), then
\(Q-X\approx Q-Y\).
\end{defn}

A Z-set of a space $X$ is a subspace $Z$ for which there is a homotopy \(H_t:X\to X\)
such that \(H_0=\id_X\) but \(H_t(X)\subset X-Z\) for all \(t>0\).
Embedding a compactum $X$ as a Z-set in $Q$ is easy: one simply embeds $X$ in
\[
	\{0\}\times\Pi_{i=2}^\infty[0,1]
\le
	\Pi_{i=1}^\infty[0,1]
=
	Q.
\]
Similarly, finite dimensional compacta can be embedded as Z-sets of finite dimensional cubes.
For a proof that a finite dimensional compactum can embedded in a finite dimensional cube, see
\cite[Th 50.5]{Mu}; the proof of the infinite dimensional case is similar.\\

It is a standard fact that homotopy equivalence implies shape equivalence
(see \cite[Ch I,Sec 4.1]{MS}).\\

\subsection{Cell-Like Equivalence}
A compactum $X$ is said to be \textit{cell-like} if it is shape
equivalent to a point.
In particular, contractible compacta are cell-like.
A map \(X\to Y\) is called \textit{cell-like} if it is surjective and
the preimage of every point is a cell-like compactum.
\footnote{Note that our definition implies that cell-like maps are proper.}\\

We say that two compacta $X$ and $Y$ are \textit{cell-like equivalent}
and write \(X\CE Y\) if there is a zig-zag of compacta and cell-like maps
\[           
\begin{array}{ccccccccccccc}
	  &          & K_1 &          &     &          & K_3 &          &     &          & K_n &          &     \\
	  & \swarrow &     & \searrow &     & \swarrow &     & \searrow & ... & \swarrow &     & \searrow &     \\
	X &          &     &          & K_2 &          &     &          &     &          &     &          & Y.
\end{array}
\]
If all compacta in this zig-zag are finite dimensional, then we say that $X$ and $Y$ are
\textit{cell-like equivalent through finite dimensions}, and write \(X\fCE Y\).\\

\subsection{The Finite Dimensional Category}
If we restrict ourselves to the category of finite-dimensional compacta, then it is known that
cell-like equivalence (that is, cell-like equivalence through finite dimensions) is strictly stronger
than shape equivalence and strictly weaker than homotopy equivalence (denoted $\HE$).  Specifically,
we have the following for finite dimensional compacta $X$ and $Y$.
\begin{facts}
\hspace{0cm}\\
\begin{enumerate}
\item $X\HE Y\then X\fCE Y$ (proven by S. Ferry in \cite[Th 2]{Fe}).\\
\item $X\fCE Y\not\then X\HE Y$.\\
\item $X\fCE Y\then X\SH Y$ (proven by R. B. Sher in \cite{Sh}).\\
\item $X\SH Y\not\then X\CE Y$ (S. Ferry gave a 1-dimensional counterexample in \cite{Fe2}).\\
\end{enumerate}
\end{facts}
A couple of notes about these facts: First of all, the theorem quoted in (1) does not explicitly
mention the finite dimensional case.  However, a careful analysis of the intermediate space $Z$
constructed in \cite{Fe} reveals that it does indeed have finite dimension if $X$ and $Y$ are
finite dimensional.\footnote{The author found formula (B) from \cite[Sec III.2]{HW} helpful in this analysis.}
The second fact is a standard
example; take $X$ to be the topologist's sine curve and $Y$ be a point $p$.  The map \(X\to Y\)
is cell-like, because $X$ has the shape of a point, but $X$ is certainly not contractible.\\

Finally, it is important to observe that (3) does \textit{not} hold if we leave the finite dimensional category,
as exhibited by J. Taylor in \cite{Ta}.  However, E. Swenson has shown in \cite{Sw} that all CAT(0)
boundaries are finite dimensional, which is why Theorem \ref{thm:celllike} is stated in the finite
dimensional category.\\

\section{Proofs of Theorems \ref{thm:celllike} and \ref{thm:SH2CE}}
\label{sec:proofs}
Theorem \ref{thm:SH2CE} follows from this next proposition together with an easy transitivity argument.

\begin{prop}
Let $X$, $Y$, and $Z$ be compacta such that \(X\SH Y\).  Then
\[
	X\ast Z\CE Y\ast Z
\]
Furthermore, if these compacta are finite dimensional, then the cell-like
equivalence may be obtained through finite dimensions.
\end{prop}

\begin{proof}
We will begin by proving the proposition without the finite dimensional hypothesis.
The proof of the finite dimensional case
is obtained by an identical argument in which $Q$ is replaced with a finite dimensional cube.\\

Imbed $X$ in $Q$ as a Z-set.  In our diagrams, we will draw $Q$ as a square with $X$ as a subsegment of
a side, as in Figure \ref{fig:XinQ}.
\begin{figure}[h]
\hspace{0.0cm}
\input{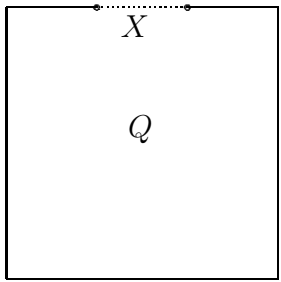}
\caption{$X\subset Q$}
\label{fig:XinQ}
\end{figure}\\
For some fixed \(z_0\in Z\), we define the space
\[
	K_3=(Q-X)\times(Z-z_0).
\]
Note that by Definition \ref{defn:shape}, we have
\[
	K_3\approx (Q-Y)\times(Z-z_0).
\]
Therefore it suffices to prove that \(X\ast Z\CE K_3^\star\), where the $\star$ denotes one-point
compactification.\\

Our cell-like equivalence zigzag between $X\ast Z$ and $K_3^\star$
will have two intermediate spaces and three cell-like maps:
\[           
\begin{array}{ccccccccccccc}
	X\ast Z  &          &     &          & K_2                       \\
	         & \stackrel{\quad\phi_1}{\searrow} &     & \stackrel{\phi_2\quad}{\swarrow} &     & \stackrel{\quad\phi_3}{\searrow}            \\
	         &          & K_1 &          &     &          & K_3^\star \\
\end{array}
\]
(see Figure \ref{fig:zigzag}).\footnote{Ric Ancel has suggested a variation on this proof which uses
only one intermediate space and two cell-like maps.}
The first intermediate space is the quotient space
\[
	K_1=X\ast Z/X\ast z_0.
\]
The other is the union of \(Q\times Z\) with a cone $\Gamma$ on the complement
of \(K_3\) (see Figure \ref{fig:conehead}).  In other words,
\begin{align*}
	K_2
=
	\Gamma\cup Q\times Z
=
	p\ast(Q\times Z-K_3)\cup Q\times Z,
\end{align*}
where $p$ denotes the cone point of $\Gamma$.
Note that \(K_3^\star=K_2/\Gamma\).  It is easy to see that $K_2$ is metrizable.
The other spaces, $K_1$ and $K_3^\star$ are metrizable because they are finite
decompositions of metrizable spaces into closed sets (see \cite{Da}).\\

The first map, $\phi_1$ is the obvious quotient map.  The last map, $\phi_3$, can also
be realized as a quotient map by writing \(K_3^\star=K_2/\Gamma\).
Both of these maps are cell-like since the only nontrivial point preimages are cones.\\

We now realize $\phi_2$ as a quotient map.
Consider the following collection of subspaces of $K_2$:
\[
	\eS=\bigl\{Q\times z\big|z\neq z_0\bigr\}\cup\bigl\{p\ast(Q\times z_0)\bigr\}
\]
and let $\phi_2$ be the quotient map onto the decomposition space \(K_2/\eS\).
Again, $\phi_2$ is obviously cell-like, since point preimages are contractible.
It suffices to prove the following claim.\\

\begin{figure}
\hspace{0cm}
\input{conehead.tex}
\caption{$K_2$}
\label{fig:conehead}
\end{figure}
\begin{figure}
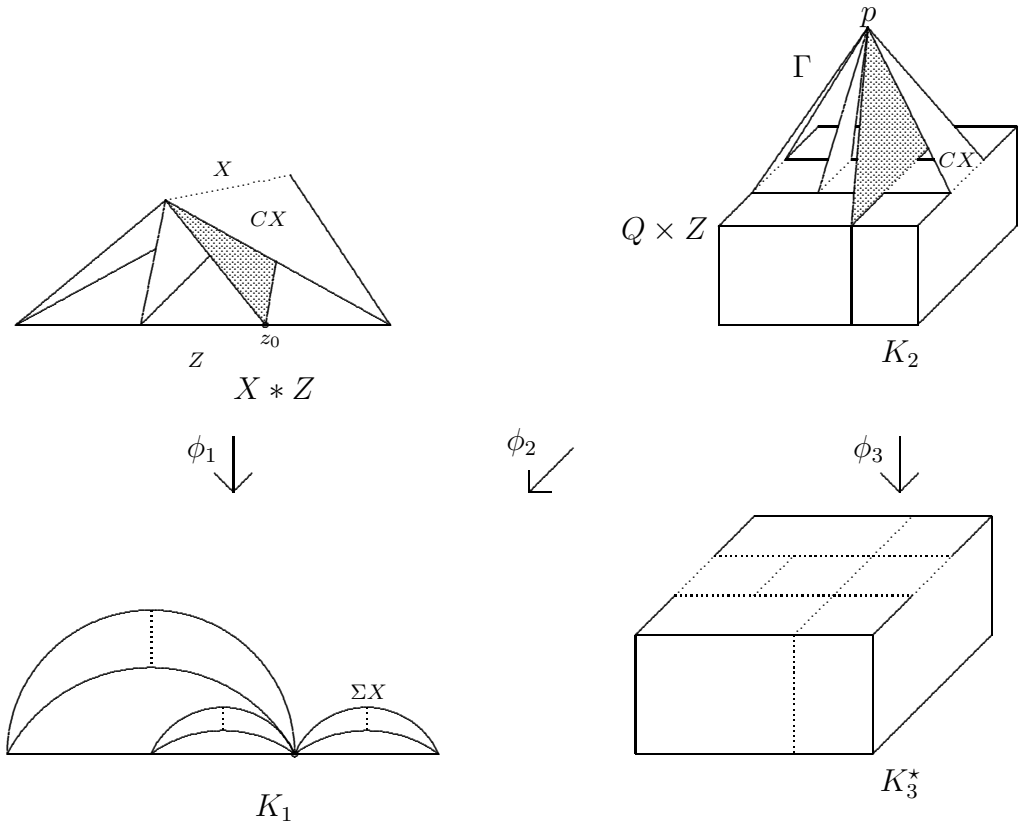

\hfill
  \begin{minipage}[t]{.45\textwidth}
    \begin{center}  
			\input{join.tex}
			$X\ast Z$
    \end{center}
  \end{minipage}
  \hfill
  \begin{minipage}[t]{.45\textwidth}
    \begin{center}  
			\input{smallconehead.tex}
			$K_2$
    \end{center}
  \end{minipage}\\
\begin{minipage}[t]{1\textwidth}
\indexspace
\begin{center}
		\hspace{0.0cm}
			\input{phi1.tex}
		\hspace{2.8cm}
			\input{phi2.tex}
		\hspace{3.2cm}
			\input{phi3.tex}\\
\end{center}
\end{minipage}\\
\hfill
  \begin{minipage}[t]{.45\textwidth}
    \begin{center}  
			\input{suspensions.tex}
			$K_1$
    \end{center}
  \end{minipage}
  \hfill
  \begin{minipage}[t]{.45\textwidth}
    \begin{center}  
			\input{deletedconehead.tex}
			$K_3^\star$
    \end{center}
  \end{minipage}
	\caption{The CE ZigZag}
	\label{fig:zigzag}
\end{figure}

\begin{claim}
	\(K_2/\eS=K_1.\)
\end{claim}

Well, \(K_2/\eS=\Gamma/\eS'\), where
\[
	\eS'=\bigl\{X\times z\big|z\neq z_0\bigr\}\cup\bigl\{p\ast(X\times z_0)\bigr\},
\]
which means that we can forget about \(Q\times Z\) for the time being.
On one hand, we can write
\[
	X\ast Z=X\times Z\times I/\sim
\]
where $\sim$ is generated by the rules
\begin{itemize}
\item [(a)] $(x,z_1,1)\sim(x,z_2,1)$ for every \(x\in X\) and \(z_1,z_2\in Z\).\\
\item [(b)] $(x_1,z,0)\sim(x_2,z,0)$ for every \(x_1,x_2\in X\) and \(z\in Z\).\\
\end{itemize}
In other words, $Z$ disappears at the top (at level 1) and $X$ disappears at the bottom (at level 0).
Passing to \(K_1\) can be thought of as adding the additional rule\\
\begin{itemize}
\item [(c)] $(x_1,z_0,t_1)\sim(x_2,z_0,t_2)$ for every \(x_1,x_2\in X\) and \(t_1,t_2\in I\),\\
\end{itemize}
which kills the cone \(X\ast z_0\).\\

On the other hand, we can write
\[
	\Gamma = X\times Z\times I/\sim
\]
where $\sim$ is generated by the rule\\
\begin{itemize}
\item [($a'$)] $(x_1,z_1,1)\sim(x_2,z_2,1)$ for every \(x_1,x_2\in X\) and \(z_1,z_2\in Z\).\\
\end{itemize}
Here \(X\times Z\) disappears at the top (at the point $p$).
Passing to \(\Gamma/\eS'\) can be thought of as adding the rules ($b$) and ($c$) from above.
But ($a$) and ($a'$) are equivalent in the context of ($c$)!
This proves the claim.\\

For the finite-dimensional version of the theorem, we simply note that
the fact that $Q$ is infinite dimensional is never used here.
Therefore we may replace $Q$ with a finite-dimensional cube and Chapman's definition with its
finite dimensional analogue proven by Chapman in \cite{Ch2} and independently by G. Venema
in \cite{Ve}.
\end{proof}

Along with Theorem \ref{thm:SH2CE}, the proof of Theorem \ref{thm:celllike} requires two other results.
The first is due to Hosaka.

\begin{Thm}
\cite[Th 2]{Ho}
Let $G=H\times K$ be a CAT(0) group with infinite factors
and $X$ be a CAT(0) $G$-space.
Then there is a CAT(0) $H$-space $Y$ and a CAT(0) $K$-space $Z$
such that
\[
	\partial X\approx\partial Y\ast\partial Z.
\]
\end{Thm}

Note that this equation is exactly what one would expect in light of the equation
\[
	\partial(Y\times Z)\approx\partial Y\ast\partial Z
\]
given in \cite[Ex II.8.11(6)]{BH}.
In fact $Y$ and $Z$ are constructed as subspaces of $X$.  The action of $H$ on $Y$ and
$K$ on $Z$ is not immediate from the original action of \(H\times K\) on $X$, however.\\

The second result is a generalization of Bestvina's theorem due to P. Ontaneda.

\begin{Thm}
\cite[Co B]{On}
Let $G$ be any CAT(0) group and $X$ and $Y$ be CAT(0) $G$-spaces.  Then \(\partial X\SH\partial Y\).
\end{Thm}

The proof of Theorem \ref{thm:celllike} is now straightforward.  Given any CAT(0) group \(G=H\times K\)
with infinite factors and any two CAT(0) $G$-spaces $X$ and $X'$, we use Hosaka's Theorem to write
\(\partial X\approx\partial Y\ast\partial Z\) and \(\partial X'\approx\partial Y'\ast\partial Z'\)
where \(Y\) and \(Y'\) are CAT(0) $H$-spaces and $Z$ and $Z'$ are CAT(0) $K$-spaces.
By Ontaneda's Theorem, we have \(\partial Y\SH\partial Y'\) and \(\partial Z\SH\partial Z'\).
Recalling that CAT(0) boundaries are always finite dimensional, we apply Theorem \ref{thm:SH2CE}
in the finite dimensional category to get that \(\partial X\fCE\partial X'\).\\

In closing, we note that the reason for requiring both factors to be infinite is because
if one of the factors, say $H$, is a finite index subgroup of $G$, then $K$ and $G$ act
geometrically on exactly the same family of CAT(0) spaces.


\begin{thebibliography}{999}

\bibitem[AGW]{AGW} F.D. Ancel, C. Guilbault, and J. Wilson,
	\emph{The Croke-Kleiner boundaries are cell-like equivalent,}
	Preprint.

\bibitem[Be]{Be} M. Bestvina,
	\emph{Local homology properties of boundaries of groups,}
	Michigan Math. J. 43 (1996), no. 1, 123-139.

\bibitem[BH]{BH} M.R. Bridson and A. Haefliger,
	\emph{Metric Spaces of Nonpositive Curvature,}
	Springer-Verlag, Berlin, 1999.

\bibitem[Bo]{Bo} K. Borsuk,
	\emph{Theory of Shape,}
	Mathematical Monographs, Vol. 59, Polish Scientific Publishing House, Warsaw. (1975).

\bibitem[BR]{BR} P.L. Bowers and K. Ruane,
	\emph{Boundaries of nonpositively curved groups of the form \(G\times\mathbb{Z}\),}
	Glasgow Math. J. 38 (1996), 177-189.

\bibitem[Ch]{Ch} T.A. Chapman,
	\emph{Lectures on Hilbert cube manifolds,}
	C.B.M.S. Regional Conf. Series in Math. 28, 1976.

\bibitem[Ch2]{Ch2} T.A. Chapman,
	\emph{Shapes of finite-dimensional compacta,}
	Fund. Math. 76 (1972), 329-353.

\bibitem[CK]{CK} C. Croke and B. Kleiner,
	\emph{Spaces with nonpositive curvature and their ideal boundaries,}
	Topology 39 (2000), no. 3, 549-556.

\bibitem[CK2]{CK2} C. Croke and B. Kleiner,
	\emph{The geodesic flow of a nonpositively curved graph manifold,}
	Geom. Funct. Anal. 12 (2002), no. 3, 479--545.

\bibitem[Da]{Da} R. J. Daverman,
	\emph{Decompositions of manifolds,}
	Pure and Applied Mathematics, 124. Academic Press, Inc., Orlando, FL, 1986. xii+317 pp.

\bibitem[Fe]{Fe} S. Ferry,
	\emph{Homotopy, simple homotopy, and compacta,}
	Topology 19 (1980), no. 2, 101--110.

\bibitem[Fe2]{Fe2} S. Ferry,
	\emph{Shape Equivalence does not Imply CE Equivalence,}
	Proc. Amer. Math. Soc. 80 (1980), no. 1, 1980 154--156.

\bibitem[Has]{Has}
	\emph{Suspensions of strong shape equivalences are CE equivalences,}
	Proc. Amer. Math. Soc. 87 (1983), no. 4, 743--745.

\bibitem[Ho]{Ho} T. Hosaka,
	\emph{On splitting theorems for CAT(0) spaces and compact geodesic spaces of non-positive curvature,}
	Preprint, arXiv:math/0405551.

\bibitem[Hr]{Hr} C. Hruska,
	\emph{Geometric invariants of spaces with isolated flats,}
	Topology 44 (2005), no. 2, 441--458. 

\bibitem[HW]{HW} W. Hurewicz and H. Wallman,
	\emph{Dimension Theory,}
	Princeton University Press, Princeton, 1969.

\bibitem[Mo]{Mo} C. Mooney,
	\emph{Examples of Non-Rigid CAT(0) Groups from the Category of Knot Groups,}
	Preprint, arXiv:math/0706.1581.

\bibitem[Mo2]{MoCK} C. Mooney,
	\emph{Generalizing the Croke-Kleiner Construction,}
	Preprint.

\bibitem[Mu]{Mu} J. R. Munkres,
	\emph{Topology, 2nd Ed,}
	Prentice-Hall, Inc., Upper Saddle River, N.J., 2000.

\bibitem[MS]{MS} S. Marde\v{s}i\'{c} and J. Segal,
	\emph{Shape Theory,}
	North-Holland, Amsterdam (1982).

\bibitem[On]{On} P. Ontaneda,
	\emph{Cocompact CAT(0) Spaces are Almost Geodesically Complete,}
	Topology 44 (2005), no. 1, 47--62.

\bibitem[Ru]{Ru} K. Ruane,
	\emph{Boundaries of CAT(0) groups of the form \(\Gamma=G\times H\),}
	Topology Appl. 92 (1999), 131-152.

\bibitem[Sh]{Sh} R. B. Sher,
	\emph{Realizing cell-like maps in Euclidean space,}
	General Topology and Appl. 2 (1972), 75--89.

\bibitem[Sw]{Sw} E. L. Swenson,
	\emph{A cut point theorem for CAT(0) groups,}
	Jour. Diff. Geom. 53 (1999) 327-358.

\bibitem[Ta]{Ta} J. Taylor,
	\emph{A counterexample in shape theory,}
	Bull. Amer. Math. Soc. 81 (1975), 629--632.

\bibitem[Ve]{Ve} G. Venema,
	\emph{Embeddings of compacta with shape dimension in the trivial range,}
	Proc. A.M.S. (1976).

\bibitem[Wi]{Wi} J.M. Wilson,
	\emph{A CAT(0) group with uncountably many distinct boundaries,}
	J. Group Theory 8 (2005), no. 2, 229--238.

\end{thebibliography}
\end{document}